# Risk-based path planning for autonomous vehicles

Qiannan Wang and Matthias Gerdts

*Abstract*— In this paper, a risk map-based path planning algorithm is introduced for autonomous vehicles. Multivariate B-splines are implemented to generate a risk map, which measures the risk of colliding with different objects. In the following step, a two-level optimal control problem is designed. At the first level, an overall lowest risk trajectory is found. Then in the second level, among the paths whose risk value is the lowest, the one with the minimum steering effort is determined. Finally, numerical simulations are carried out with the optimal control software OCPID-DAE1. Results show that this method is useful and meaningful to find a path of minimum risk for an autonomous vehicle.

## I. Introduction

Nowadays, autonomous driving is an attractive and hot topic for both manufacturers and researchers. It is treated as a new trend and challenge for automobile industry, which will combine interdisciplinary technology together. Autonomous driving is divided into several levels, from "no automation" to "full automation". After several decades' efforts, amounts of advanced driver assistance systems (ADAS) have been designed for assisted and partial automation, e.g. emergency braking, adaptive cruise control and self-parking systems. In order to achieve fully automated driving, further concentration should be focused on urban area, where environment is more complicated and safety problems are more urgent.

Path planning is a necessary step for autonomous driving. It works like "brain" of the vehicle by defining the way to destination, while satisfying constraints at the same time. It has been widely studied and a large number of methods have been developed. Some popular methods are: (1) graph search based, like Dijkstra's algorithm introduced in [1-2], A* method developed in [3]; (2) sampling based, e.g. Probabilistic Roadmap Method (PRM) [4] and Rapidly-exploring Random Tree (RRT) [5-6]; (3) interpolating curve based, e.g. spline curves [7], polynomial curves [8-9]; (4) numerical optimization methods, e.g. model predictive control (MPC) [10-11]). Generally, graph search based and sampling based methods are suitable for finding a global path. Paths generated by interpolating curve-based methods are usually not continuous. Optimization-based methods are available for local path planning and are capable of taking into account constraints. The problem to be solved in this paper is related with local path planning. Thus, a gradient-based numerical optimization method is chosen.

*This Research is funded by funded by dtec.bw - Digitalization and Technology Research Center of the Bundeswehr [EMERGENCY-VRD], https://dtecbw.de/.

Qiannan Wang is research assistant in Institute of Applied Mathematics and Scientific Computing, Universität der Bundeswehr München, Germany (e-mail: qiannan.wang@unibw.de).

Matthias Gerdts is the professor of Institute of Applied Mathematics and Scientific Computing, Universität der Bundeswehr München, Germany (e-mail: matthias.gerdts@unibw.de).

When autonomous vehicles move in uncertain or populated environments, risk should be accounted for to produce an effective and safe path. Risk-aware based path planning has been investigated by some researchers for unmanned aerial vehicles [12-15]. The general idea is to compute an effective path that minimizes the risk to the population. In [12], path planning is performed in two phases: offline based on static risk factors and online based on dynamic information. A* method is used for offline and an algorithm called Borderland is developed to adjust the portion of path. In [13], a risk assessment technology and bi-objective optimization methods are designed to find a low-risk and less time solution. These papers deal with the problem of finding a lowest risk path, where the risk is only distinguished between 0 and 1. That is, if population appears, the risk value is 1. Otherwise, it is 0. However, this method of modeling risk map is not suitable for autonomous vehicles. Normally autonomous vehicles have to face various objects in complicated environment, e.g. pedestrians, bicyclists, trees, and other vehicles. Assigning the same risk value for all kinds of participants is not very helpful. To generate a collision-free or decrease the severity of collisions, disparate risk values for different objects should be well defined. Therefore, a new method of modelling risk map and optimization for path planning will be discussed here.

In order to achieve fully autonomous driving in urban environment, autonomous vehicles should have the ability of avoiding collisions - if possible at all - with any other objects in the area. When collisions cannot be avoided, it has to release the harm to the greatest extent. To realize this goal, a risk map based path planning method is investigated. The main contributions of this work are:

(1) A risk map is generated based on multivariate B-splines. Here, all appeared objects are classified and corresponding risk values are precisely defined. Through a least-squared problem, the risk map is generated.

(2) A two-level optimal control problem based on the risk map is designed for path planning, where the lowest risk and the minimum steering efforts are under consideration.

The structure of the paper is constituted as follows. The method of generating risk map is mentioned in section II. Following in section III, detailed information of the path planning method is presented. Simulation results are discussed in section IV and conclusions and future work are given in section V.

## II. Generating risk map

In order to generate a lowest risk path, a well-defined risk map for each scenario is essential. The basic idea to generate a risk map is that, by giving risk values of some points in the scenario, the risk values of the whole area are acquired. To this end the various objects in the two-dimensional traffic scenario will be assigned risk values risk(x,y) depending on their positions (x,y). These objects are, e.g., the lanes of the roads, other vehicles, pedestrians, trees, shoulder. Technically these objects are detected by suitable sensors such as radar, LIDAR, or cameras. These sensors deliver scattered data points which are labelled. The label indicates the detected object and size and shape can be estimated from the data. The first step in the derivation of a risk map risk(x,y) is to generate a sufficiently smooth function that approximates the scattered data. A smooth approximation is required by the gradient-based optimization methods to be used for path planning.

*Remarks: In this paper for simplicity we restrict the discussion to risk maps which solely depend on the position (x,y) of the objects in the traffic scenario. However, it is straightforward, but technically more involved, to consider risk maps which take into account further quantities of influence, e.g., collision angle and geometry of objects. Likewise we only consider static risk maps. Again, adding the time as additional variable in the risk map is straightforward.*

Given scattered data points, the first task is to generate a smooth risk map. This could be achieved by polynomial interpolation or spline interpolation. Owing to the large number of data points and the presence of perturbations, this would lead to oscillations and thus is unsuitable for our purposes. Instead we use multivariate B-spline approximations for generating the risk map. These allow to control the smoothness of the approximations conveniently by choosing the orders of the elementary B-splines. In addition, Runge's phenomenon can be avoided, in which oscillations can occur between points if high-degree polynomials are used.

### A. Multivariate B-spline Approximation

Let M scattered data points

$$(x_s, y_s, r_s)^T \in \mathbb{R}^3, s = 1 \dots M, \qquad (1)$$

be given. M is a potentially large number. It is assumed that these data points are realizations of a function $r : \mathbb{R}^2 \to \mathbb{R}$ with $r_s = r(x_s, y_s)$, $s = 1 \dots M$. We intend to find an approximation *risk* of $r$ with $risk(x_s, y_s) \approx r_s$, $s = 1, \dots, M$. To this end we express *risk* as a linear combination of elementary B-spline basic functions of orders $k_1$ and $k_2$, that is

$$risk(x,y) := \sum_{i=1}^{n_1} \sum_{j=1}^{n_2} c_{ij} B_i^{k_1}(x) B_j^{k_2}(y) \qquad (2)$$

Where, $i$ and $j$ are the number of knots in *x* and *y* directions. The coefficients $c_{ij}$ ($i=1\dots n_1, j=1\dots n_2$) are to be determined by solving the least-squares problem in (3) with respect to $c \in \mathbb{R}^{n_1 \times n_2}$.

$$\text{Minimize } \tfrac{1}{2} \sum_{l=1}^{M} (r_l - risk(x_l, y_l))^2 + \tfrac{\lambda}{2} \|c\|^2 \qquad (3)$$

Where $\lambda > 0$ is a regularization parameter. The elementary basic functions $B_i^k$ of order $k \in N$ are defined as follows. For given data points $t_0 < t_1 < \cdots < t_N$, we introduce an auxiliary grid

$$G_N^k := \{\tau_i | i = 1, \dots, N + 2k - 1\} \qquad (4)$$

with grid points

$$\tau_i := \begin{cases} t_0, & \text{if } 1 \le i \le k \\ t_{i-k}, & \text{if } k+1 \le i \le N+k-1 \\ t_N, & \text{if } N+K \le i \le N+2k-1 \end{cases} \qquad (5)$$

The *i*-th B-spline of order *k* is defined recursively by

$$B_i^1(t) := \begin{cases} 1, & \text{if } \tau_i \le t < \tau_{i+1} \\ 0, & \text{otherwise} \end{cases}, \qquad (6)$$

$$B_i^k(t) := \tfrac{t-\tau_i}{\tau_{i+k-1}-\tau_i} B_i^{k-1}(t) + \tfrac{\tau_{i+k}-t}{\tau_{i+k}-\tau_{i+1}} B_{i+1}^{k-1}(t) \qquad (7)$$

for $i=1,\dots, N+k-1$. B-splines have a local support given by the interval ($\tau_i$, $\tau_{i+k}$) with

$$B_i^k(t) \begin{cases} > 0, & \text{if } t \in (\tau_i, \tau_{i+k}), \\ = 0, & \text{otherwise}, \end{cases} \text{ for } k>1. \qquad (8)$$

To generate risk map, we only need non-negative values. Therefore, we consider the constrained least-squares problem:

$$\text{Minimize } \quad J(c) := \tfrac{1}{2} \|Ac - z\|^2 + \tfrac{\lambda}{2} \|c\|^2 \qquad (9)$$

with c ≥ 0. This is a convex problem. If $\lambda > 0$, it is even strictly convex. The necessary and sufficient optimality conditions for a solution *c* of this constrained optimization problem read as follows:

$$c = Proj_{\mathbb{R}_+^n}(c - \gamma((A^T A + \lambda I)c - A^T z)) \qquad (10)$$

where γ >0 is an arbitrary number and the projection is onto the set of non-negative vectors of appropriate dimension. This non-smooth equation can be solved with the fixed-point iteration

$$c^{(k+1)} = Proj_{\mathbb{R}_+^n}(c^{(k)} - \gamma((A^T A + \lambda I)c^{(k)} - A^T z)) \qquad (11)$$

where k= 0,1,2…. To this end, $\gamma \in (0,1)$ has to be chosen appropriately.

### B. Risk values of different objects

In subsection A, we have discussed how to use multivariate B-splines to approximate the risk function. In this work, it is assumed that information of objects such as type of participants and location are already known through environment perception system. It will not be discussed here. To generate a reasonable risk map, it is important to give suitable risk values for objects that appear in the scenario. It has to be pointed out that the actual risk values determine the risk-minimal path and thus the actual behavior of the autonomous vehicle. Hence, assigning the risk values is a crucial task, especially in situations where collisions are unavoidable. Precise risk values should be evaluated by psychologists based on ethical decisions.

In this paper, we take intersections as an example. Intersection scenarios are typically complex, since there may be involved various participants such as pedestrians, cars, trucks, and trees. Furthermore, road shoulders should also be considered. At the moment, there are no exact rules for definition of risk values. For illustration and simulation, we propose a rule as shown in Table I. The general idea is that, for pedestrians, the value is the highest. Trucks have relative lower values than pedestrians but still higher than cars. Trees and road shoulders have the lowest values. As pointed out before, these values need a thorough discussion and evaluation before they are implemented in a real system.

Moreover, the shape of the objects should also be taken into account, including the shape of the ego vehicle. By giving different risk values for each object, the 'behavior' or 'decision' of the autonomous vehicle will be influenced.

For the first step investigation, only static obstacles are considered. The extension to dynamic obstacles is also possible with some technical efforts regarding the multivariate B-splines. Based on the description in sections A and B, a risk map can be generated. Figure 1 is a schematic diagram of a risk map for an intersection, where a truck, a pedestrian and a passenger car are present.

TABLE I.   RISK VALUES OF VARIOUS OBJECTS

| Objects | Risk value |
|---|---|
| pedestrians | 40 |
| trucks | 30 |
| cars | 20 |
| trees | 10 |
| road shoulders | 10 |

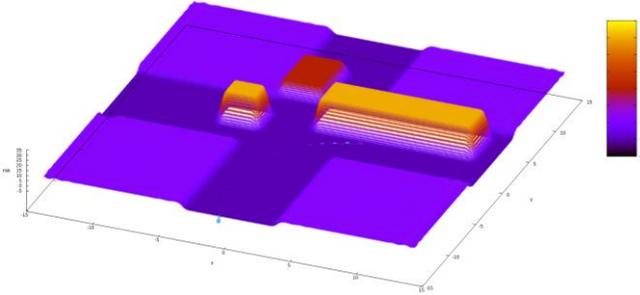

Figure 1.   Example of a risk map at an intersection scenario

### III. PATH PLANNING CONTROL

After getting the risk map, now we go to the core point, designing the path planning problem. In this section, the vehicle model and the optimal control problem are introduced. We like to emphasize that we are particularly interested in crucial scenarios where collisions are likely or cannot be avoided at all. The duration of such scenarios is typically short in the range of a few seconds, while the velocities are comparatively high such that pure braking is not feasible.

#### A. Vehicle model

A kinematic vehicle model is used in this work. The equations of motions are shown as the following:

$$\dot{x} = v\cos\varphi \quad (12)$$

$$\dot{y} = v\sin\varphi \quad (13)$$

$$\dot{\varphi} = v\tan\delta/L \quad (14)$$

$$\dot{v} = a \quad (15)$$

$$\dot{\delta} = (\delta_s - \delta)/\Delta T \quad (16)$$

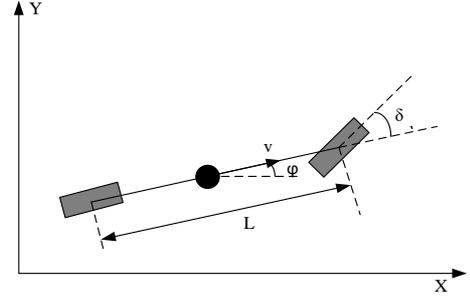

Figure 2.   Vehicle kinematic model

Here, *x* and *y* are coordinates of the car´s reference point (midpoint of rear axle) and *v* denotes the velocity of the vehicle. The yaw angle is φ and the steering angle is *δ*. *L* represents the distance from the front axle to rear axle of the vehicle and the control input *a* is the acceleration of the vehicle. $\delta_s$ is the desired steering angle, which is a control input. The constant value $\Delta T$ is used to model a latency in the steering device.

#### B. Optimal control problem

In section II, the risk function has been analyzed. We use it to formulate optimal control problems for path planning.

First we solve the following optimal control problem ($OCP_1$), which aims at minimizing the total risk. Note that the risk map is weighted with the velocity since the severity of impacts increases with the velocity:

Minimize

$$J_1 = \int_{t_0}^{t_f}(v(t) * \dot{risk}(x(t),y(t)))^2 dt \quad (17)$$

subject to system equations (12) to (16), control constraints

$$\begin{cases} -10 \leq a \leq 2 \\ -0.4 \leq \delta_s \leq 0.4 \end{cases} \quad (18)$$

and initial conditions

$$(x(0),y(0),\varphi(0),v(0),\delta(0)) = (x_0,y_0,\varphi_0,v_0,\delta_0). \quad (19)$$

Here, $t_0$ and $t_f$ denote the initial and final time, respectively. Through the minimization of the risk measure $J_1$, we obtain a risk minimal path. Let $J_1^*$ denote the optimal total risk value of ($OCP_1$). In general the solution is not unique and several paths with minimal risk measure exist. Among all such risk minimal paths we try to identify a suitable path for the autonomous vehicle.

To this end, in a second step we solve the following optimal problem ($OCP_2$), which aims at finding a trajectory with minimized steering efforts among all trajectories with the same minimum total risk:

Minimize

$$J_2 = \int_{t_0}^{t_f} \delta_s(t)^2 dt \quad (20)$$

subject to system equations (12) to (16), constraints for acceleration and steering angle in (18), the initial conditions (19), and the additional risk constraint

$$\int_{t_0}^{t_f}(v(t) * \dot{risk}(x(t),y(t)))^2 dt \leq J_1^* + \varepsilon. \quad (21)$$

Herein, $\varepsilon \geq 0$ is a small value, which can be used to slightly relax the optimal risk measure $J_1^*$ from $OCP_1$. This relaxation might increase numerical stability.

## IV. SIMULATION AND RESULTS

In this section, the designed path planning method will be implemented in simulation software OCPID-DAE1 (Optimal Control and Parameter Identification with Differential Algebraic Equations of Index 1) for validation. OCPID-DAE1 is an optimal control tool, which implements a direct shooting method. The resulting nonlinear optimization problem is being solved by a sequential quadratic programming (SQP) method. Detailed information can be found in user manual in [16]. The intersection scenario is chosen for investigation, due to its complicated environment and high potential for accidents. At intersections, private cars, trucks, pedestrians and trees along the road boundaries are considered. Critical scenarios, where it is congested or collisions cannot be avoided, are taken into account. At this stage, all obstacles are static for simplicity, but an extension towards non-static obstacles is straightforward with some technical effort. The initial states of the ego vehicle in all cases are the same, $(x(0),y(0),\varphi(0),v(0),\delta(0))=(27, 17, \pi, 20, 0)$.

### A. Case 1

In section III, a two-level optimal problem is proposed. In order to validate the effectiveness of this design, a comparison with one level optimal control problem is carried out. The one level control problem only considers the minimized risk value, and the objective function 1 is implemented. The scenario is shown in Figure 3. The ego vehicle plans to traverses the crossing. There are two queues of pedestrians, one truck, one private car and some trees. The ego vehicle has to find a path among all these objects.

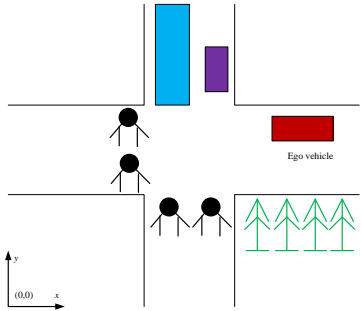

Figure 3. Scenario in case 1

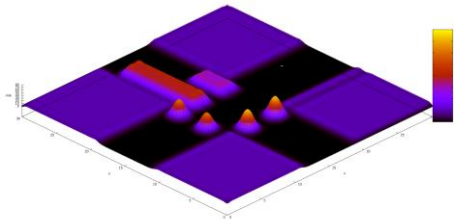

Figure 4. Risk map for Case 1

The risk map for case 1 is indicated in Figure 4. Here, pedestrians are modeled as circles, truck and car are modeled as rectangles. The paths generated by these two controllers are shown in Figure 5 and 6, respectively. From the results, it seems that there are not so many differences.

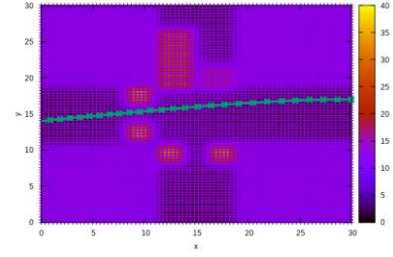

Figure 5. Path generated by two-level controller

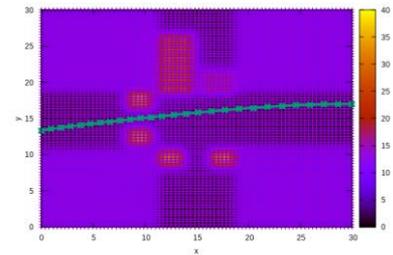

Figure 6. Path generated by one level controller

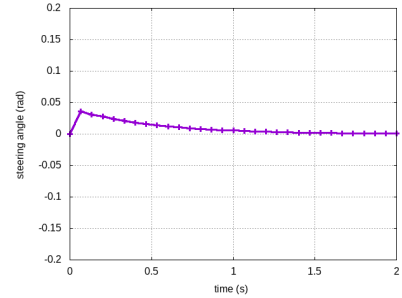

Figure 7. Steering angle by two-level controller

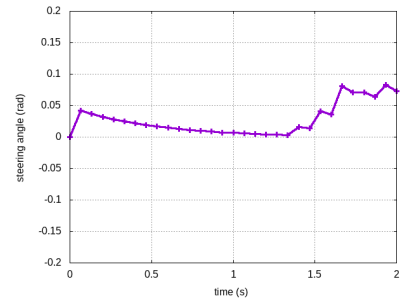

Figure 8. Steering angle by one level controller

Steering angles, accelerations, velocities and risk values of these two controllers are show in Figure 7 to 14. Change of deceleration and velocity of these two controllers keep almost the same. The steering angle of the two-level controller is better than one level controller, due to a further control on steering. The risk value of the two-level controller is a litter higher than that of the one level controller. It is because that there is a constant adjustment parameter for risk value constraints in the two-level controller. As a conclusion, both one level and two-level controller can generate feasible paths. The two-level controller performs better by minizing steering effort. Later in case 2 and 3, only the two-level controller is used.

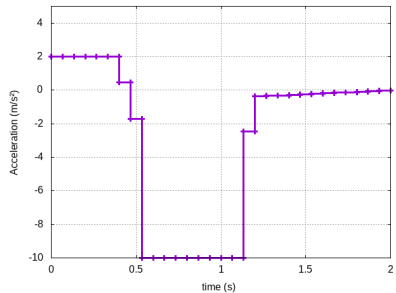

Figure 9. Acceleration by two-level controller

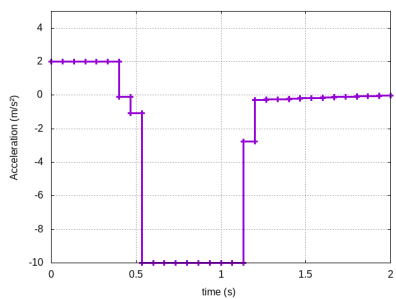

Figure 10. Acceleration by one level controller

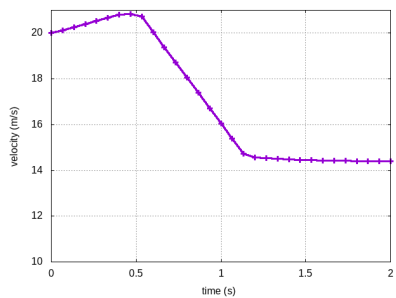

Figure 11. Velocity by two-level controller

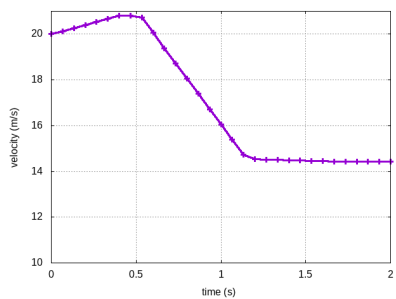

Figure 12. Velocity by one level controller

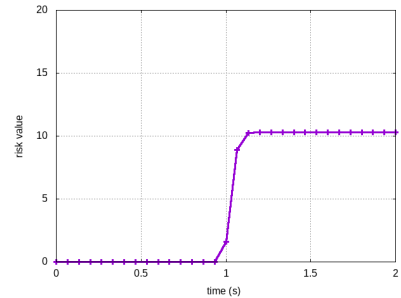

Figure 13. Risk value by two-level controller

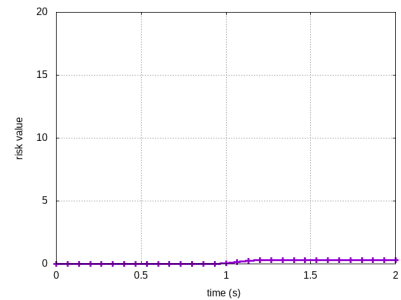

Figure 14. Risk value by one level controller

### B. Case 2

In case 2, the scenario is the same as shown in Figure 3. We change the shape of pedestrians. The two queues of pedestrians are modeled as two rectangles as shown in Figure 15.

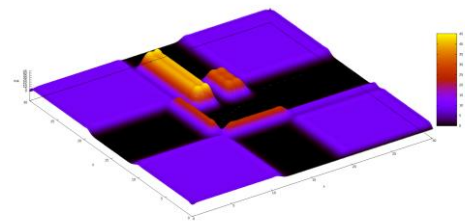

Figure 15. Risk map of case 2

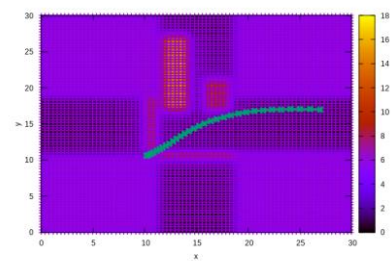

Figure 16. Generated path in Case 2

From the generated path in Figure 16, it can be seen that, the ego vehicle goes through the gap between two queues of pedestrians, where the risk is relative lower. Compared with the path generated in Figure 5, the result is totally different as a result of varied modelling shape. The velocity is shown in Figure 17. It decreases all the way to finally pass through. The deceleration shows that the vehicle brakes all the time. The steering angle changes as shown in Figure 19. The risk value is shown in Figure 20. It shows that there are no free spaces and the vehicle tries to go along a lowest risk path. From the simulation results, it is concluded that, with the help of risk map, the ego vehicle can find a relative lowest risk value path.

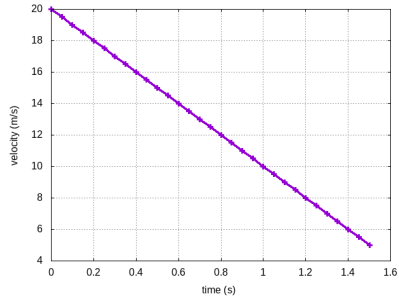

Figure 17. Velocity in Case 2

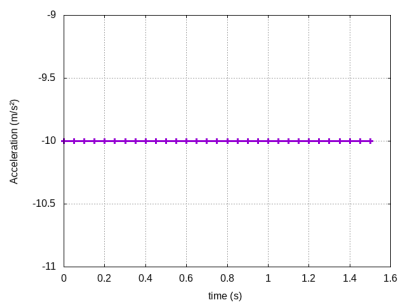

Figure 18. Deceleration in Case 2

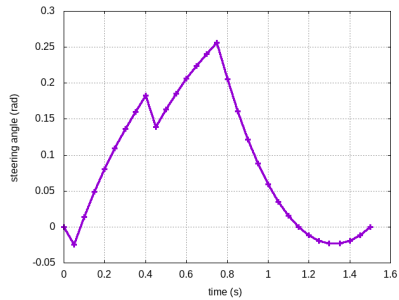

Figure 19. Steering angle in Case 2

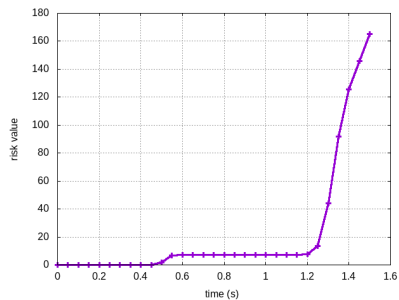

Figure 20. Risk value in case 2

*C. Case 3*

In this case, a more complicated scenario is studied. Based on the objects appear in case 1, another private car occurs, which makes the intersection more crowded. It will be more challengeable for the ego vehicle to find a path.

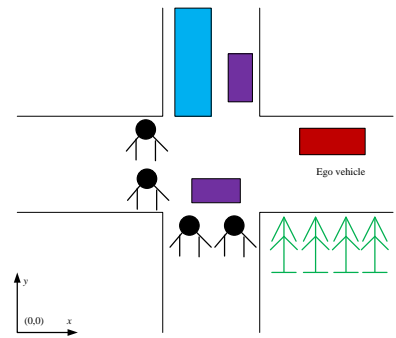

Figure 21. Case 3 Scenario

The risk map is shown as in Figure 22. The generated path is shown in Figure 23 (green line). From the path, it can be seen that, the ego vehicle goes through the two pedestrian. As the shapes of objects are already considered, this path is collision-free.

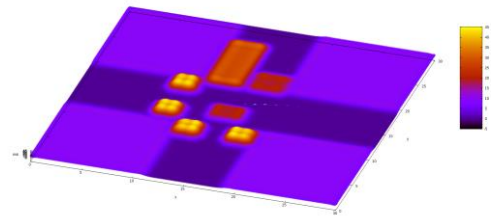

Figure 22. Risk map of Case 3

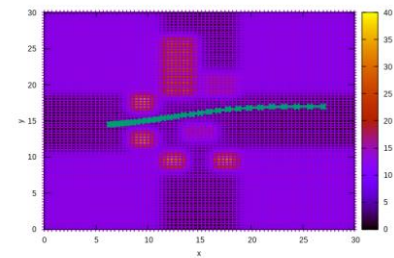

Figure 23. Generated path in case 3

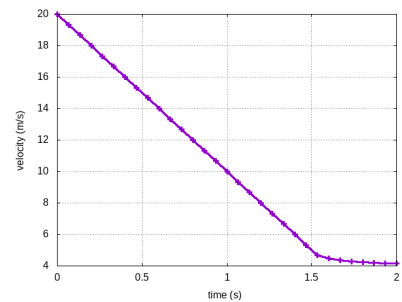

Figure 24. Velocity in case 3

The velocity is shown in Figure 24. It can be observed that it decreases all the way to finally pass through. After it passes, the velocity keeps almost constant. The deceleration shows that the vehicle firstly brakes at the maximum value. At round 1.5s, it begins to increase to 0. The steering angle is shown in

Figure 26. The risk value is shown in Figure 27. The risk value keeps around 0, which means no collision with any objects.

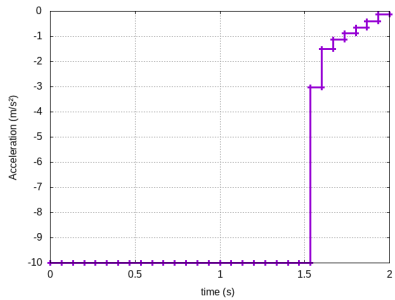

Figure 25. Deceleration in case 3

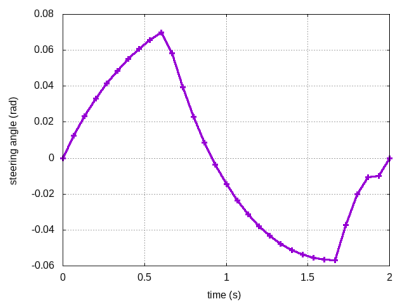

Figure 26. Steering angle in case 3

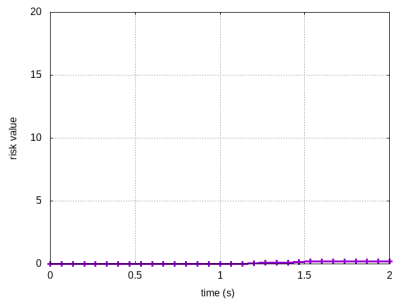

Figure 27. Risk values in case 3

## V. CONCLUSIONS AND FUTURE WORK

In this paper, a risk map based path planning method is designed. Evaluation of risk values of different objects are studied, which are important for decision making, especially for some critical situations. Risk value and steering efforts are minimized to get an optimal path. Simulation results show that the method works well and it is a promising way to deal with complicated autonomous driving scenarios in the future.

For future work, further development of risk map with dynamic obstacles should be in the focus. Investigation of influence of different risk values on the performance of the vehicle is an interesting topic. Moreover, implementation of this algorithm on a real car is also a perspective.

### ACKNOWLEDGMENT

This research is funded by detc.bw- Digitalization and Technology Research Center of the Bundeswehr München. This work is part of the project [EMERGENCY-VRD], detailed information can be found at https://dtecbw.de/.